\documentclass[a4paper,11pt,parskip=half]{amsart}
\usepackage[utf8]{inputenc}

\usepackage{amsmath, amsthm, amssymb}
\usepackage{graphicx}
\usepackage{subcaption}
\usepackage{upgreek}
\usepackage{mathrsfs}
\usepackage{xcolor}
\usepackage{mathtools}
\usepackage[foot]{amsaddr}
\usepackage[T1]{fontenc}
\usepackage{palatino}
\usepackage[margin=1 in]{geometry}
\usepackage[colorlinks=true]{hyperref}
\hypersetup{urlcolor=blue, citecolor=red}

\DeclarePairedDelimiter{\abs}{\lvert}{\rvert}
\DeclarePairedDelimiter{\set}{\{}{\}}

\DeclarePairedDelimiter{\prt}{(}{)}
\DeclarePairedDelimiter{\brk}{[}{]}

\DeclareMathOperator{\diver}{div}

\newcommand{\R}{{\mathbb R}}

\DeclareMathAlphabet{\mathup}{OT1}{\familydefault}{m}{n}
\newcommand{\dx}[1]{\mathop{}\!\mathup{d} #1}

\newcommand{\ddt}{\frac{\dx{}}{\dx{t}}}
\newcommand{\partialt}[1]{\frac{\partial #1}{\partial t}}

\newcommand{\fpartial}[1]{\frac{\partial}{\partial #1}}
\newcommand{\partialxx}[1]{\frac{\partial^2}{\partial x^2}}

\newcommand{\wm}{\abs{w}_{-}}

\newcommand {\f}  {\frac}

\newcommand {\al} {\alpha}

\newcommand{\ie}{\emph{i.e.}}
\newcommand{\cf}{\emph{cf.}\;}
\newcommand{\e}{\varepsilon}

\theoremstyle{plain}
\newtheorem{theorem}{Theorem}[section]
\newtheorem{lemma}[theorem]{Lemma}

\theoremstyle{remark}
\newtheorem{remark}[theorem]{\bf Remark}

\newcommand{\curlyL}{\mathcal{L}}
\newcommand{\qqand}{\qquad \text{and} \qquad}
\newcommand{\qand}{\quad \text{and} \quad}

\begin{document}

\title{The Aronson-B\'enilan Estimate in Lebesgue spaces}
\date{\today}

\author{Giulia Bevilacqua$^{1}$}
\author{Beno\^it Perthame$^{2}$}
\author{Markus Schmidtchen$^{3}$}

\address{$^{1}$ MOX -- Politecnico di Milano, piazza Leonardo da Vinci 32, 20133, Milano, Italy (giulia.bevilacqua@polimi.it)}
\address{$^{2}$ Sorbonne Universit\'e, CNRS, Universit\'e de Paris, Inria, Laboratoire Jacques-Louis Lions, 75005 Paris, France,  (benoit.perthame@sorbonne-universite.fr)}

\address{$^{3}$ Sorbonne Universit\'e, CNRS, Universit\'e de Paris, Laboratoire Jacques-Louis Lions,75005 Paris, France (markus.schmidtchen@upmc.fr)}
\maketitle

\begin{abstract}
     In a celebrated three-pages long paper in 1979, Aronson and B\'enilan obtained a remarkable estimate on second order derivatives for the solution of the porous media equation. Since its publication, the theory of porous medium flow has expanded relentlessly with applications including thermodynamics, gas flow, ground water flow as well as ecological population dynamics. The purpose of this paper is to clarify the use of recent extensions of the Aronson and B\'enilan estimate in $L
    ^p$ spaces, of some modifications and improvements, as well as to show certain limitations of their strategy.
\end{abstract}

\noindent{\makebox[1in]\hrulefill}\newline
\section{Introduction}
The porous medium equation (PME) is, without a doubt, one of mathematics' evergreens as it typically occurs whenever a quantity, $n(t,x)$, evolves according to a continuity equation,
\begin{align}
    \label{eq:continuity_eqn}
    \partialt n + \diver \prt*{n u} = 0,
\end{align}
where the velocity is given by \emph{Darcy's law}, \cf \cite{Dar1856}, \ie,
\begin{align*}
    u =  -\nabla p.
\end{align*}
Here, $p$ denotes the \emph{pressure}. Since the equation, in its current form, is not closed, a \emph{constitutive pressure law}, also referred to as \emph{equation of state}, is chosen to close the equation. Such a law relates the pressure directly to the quantity, $n$, affected by the pressure, \ie, $p=p(n)$. As mentioned before, its practical applications are manifold, reaching from problems related to ground water flow \cite{Dar1856}, nonlinear heat transfer \cite{ZR66}, and population dynamics, \cite{GN75, Ito52, Mor50}, to name just a few. For an extensive and elaborate treatise of the porous medium equation, we refer the reader to the homonymous book by Vazquez \cite{Vaz07}, and references therein.

One of the most fascinating properties of nonlinear diffusion phenomena is the finite speed of propagation. While solutions of the linear diffusion equation
have an instantaneous regularising effect, \ie, solutions become positive and smooth after an arbitrarily short time, solutions of the porous medium equation exhibit a behaviour quite different from that of the linear case -- solutions remain with limited regularity and compactly supported if they were compactly supported initially, a  phenomenon often referred to as \emph{finite speed of propagation}. In \cite{OKYL58}, the authors introduce a notion of weak solutions and give an existence and uniqueness result of weak solutions to the filtration equation,
\begin{align}
    \label{eq:FiltrationEqn}
    \partialt n = \frac{\partial^2}{\partial x^2} \phi(n),
\end{align}
for a certain class of functions $\phi$.  Moreover, they show that the equation is satisfied in the classical sense in neighbourhoods of points, $(t,x)$, where $n(t,x)>0$, and, that solutions, emerging from compactly supported initial data, have compact support for all times.

In a later paper, \cite{Kal66}, more properties of the porous medium equation were shown. In particular, invaded regions will remain covered with the density, $n$, \cf \cite[Lemma 2]{Kal66} for all times, every point in space will be invaded by the density after a sufficiently long time, \cite[Lemma 3]{Kal66} and regions of vacuum do not fill up spontaneously, \cf \cite[Lemma 3]{Kal66}.

Intrigued by the fact that the support of any solution emanating from compactly supported initial data is bounded by two free boundaries, or \emph{interfaces}, \cf  \cite{Kal66}, Aronson proved a characterisation of the free boundary speed which is directly related to the pressure gradient which acts as the formal velocity in Eq. \eqref{eq:continuity_eqn}, see \cite{Aro70a}. In this paper, it is assumed that some initial data supported on an interval $I=(a_1, a_2)$ are given.  In order to establish the aforementioned characterisation of the speed of the moving boundary, Aronson remarks that some control of the quantity ``$p_{xx}$'' is required for the analysis. However, the author was able 
to construct a counterexample explaining that this type of regularity cannot, in general, be expected, which was already known in the case of the explicit Barenblatt-Pattle solution, discovered in 1952, \cf \cite{Bar52, Pat59}.  In fact, in \cite{Aro70}, he provides smooth initial data (of $C^\infty$-regularity) that exhibit  blow-up of ``$p_{xx}$'' in finite time. Assuming a Power-law for the pressure, \ie , 
\begin{align}
    p = n^{\gamma},
\end{align}
Aronson investigates the behaviour of the following one dimension Cauchy problem 
\begin{equation}
\label{eq:counterexample}
\left\{
\begin{aligned}
&\frac{\partial p}{\partial t}  = \gamma p p_{xx} +  p_{x}^2\\
&p(x,0) = \cos^2(x) 
\end{aligned}
\right.
\end{equation}
and he proves that, in general, it is not possible to estimate the second derivatives of the solution of Eq. \eqref{eq:counterexample}
in terms of the bounds for the derivatives of the initial data. Indeed, he obtains that the second derivative of the pressure behaves like
\begin{align}
    p_{xx} = \frac{2 T}{T-t}, \qquad \hbox{where} \qquad T = \frac{m-1}{2m (m-1)},
\end{align}
which implies that it exceeds any bound in finite time, for more details one can refer to \cite[Theorem 1]{Aro70}. Moreover, a similar result can be obtained when the initial data has compact support, \cf  \cite[p. 301]{Aro70}, Example 2.

Thus a different type of control is necessary. In \cite[Lemma 2]{Aro70a}, Aronson proves that if
\begin{align}
    \mathrm{ess\, inf}_{I} \frac{\partial^2 p}{\partial x^2}(0,x) \geq - \alpha,
\end{align}
for some $\alpha \geq 0$, then 
\begin{align}
    \frac{\partial^2 p}{\partial x^2}(t,x) \geq - \alpha,
\end{align}
for all $(t,x)$ such that $n(t,x)> 0$. To the best of our knowledge, this is the first time this type of lower bound on the Laplacian\footnote{At this time, the results by Oleinik, Kalashnikov, and Aronson only concern the one-dimensional case.} of the pressure is obtained. At the same time, this observation acts as the foundation of the more refined version with less restrictions on the initial data, obtained in \cite{AB79}, in 1979.

The lower bound on the Laplacian of the pressure is achieved by regularising the initial data in the following way:
\begin{align}
    p_0^n(x) = (k_n\star p_0)(x) + 2^{2-n} K,
\end{align}
where $K$ is the Lipschitz constant of $p_0$ and $k_n$ is a sequence of smoothing kernels converging to a Dirac delta. For the smoothened and strictly positive initial data classical solutions exist, \cf \cite{OKYL58}, and an equation for the second derivative of the pressure, $p_{xx}$, can be found. Aronson observes that the resulting equation for $p_{xx}$ can be cast into a form whose parabolic operator satisfies the maximum principle presented in \cite[Theorem 8]{IKO62}. Ultimately, this allows to deduce the uniform bound from below on the second derivative of the pressure.\\

Later, Aronson and B\'enilan show that a similar estimate (known as Aronson-B\'enilan estimate, \ie, AB-estimate) can be obtained in the multi-dimensional case, \cf \cite{AB79}. Under no additional assumptions\footnote{only exponents ``large enough", cf. their paper} they show that
\begin{align}
\label{eq:estimate_deltaP}
    \Delta p \geq - \frac{c}{t},
\end{align}
for some constant $c>0$. Let us note that the same result is already mentioned in \cite{Aro70a} as a note, since the regularising effect was not the main focus in the derivation of the boundary speed characterisation in one dimension.\\

In 1982, Crandall and Pierre generalise the Aronson-B\'enilan estimate for the initial-value problem associated to the filtration equation, \cf Eq. \eqref{eq:FiltrationEqn}, given by
\begin{equation}
    \label{eq:pmeq_phi}
    \left\{
    \begin{aligned}
    \partial_t n &= \Delta \phi(n),\\
    n(0,x) &=n_0(x),
    \end{aligned}
    \right.
\end{equation}
for $t >0$, and $x \in \mathbb{R}^N$, where $\phi$ is a non-negative, non-decreasing, continuous function with $\phi(0) = 0$, \cf \cite{CP82}.
They prove that if $\phi$ satisfies an inequality, \cf  \cite[Eq. (3)]{CP82} which formally controls the growth of $\phi$,  the solution $n$ of Eq. \eqref{eq:pmeq_phi} satisfies
\begin{align}
    \partial_t  \phi \left(n\right) \geq \frac{K}{t}\left(\phi(n) + a\right),
\end{align}
for some constants $K>0$ and $a\geq 0$.  In the Power-law case, $\phi(n) = n^\gamma$, with $a = 0$ and
\begin{align*}
    K \geq \frac{\gamma}{\gamma-1+2/N},
\end{align*}
obtained from the aforementioned inequality, \cf \cite[Eq. (3)]{CP82}, the solution $n$ underlies the same regularising effect as that of \cite{AB79}. With this approach, they are able to extend the AB-type estimate to cover a larger class of problems, including, for instance, the Stefan problem (see \cite{CE83}) which holds for $a = 1$, $K = N/2$. In their proof, they use a similar technique adopted in this work: they study the evolution of the Laplacian of the pressure, $\Delta p$, which satisfies an appropriate parabolic inequality. Upon introducing a suitable function $h(p)$ characterised a posteriori, they  deduce the time estimate on the solution $n$. Here, different from \cite{CP82}, first we extend the AB-estimate in all the interesting Lesbegue spaces and then, in $L^{\infty}$, we are able to weaken the condition on the the function $\phi$, \ie, we do not need to impose their Eq. (3). Precisely, our needed condition to get the AB-type estimate are the same as the ones in \cite{CP82}. There are two main differences: first they are not addressing the incompressible limit and then, to assume their Eq. (3), they are forced to select a specific form of the weight which also ensures that a specific power of $\phi$ is convex. In our paper, first we do not need any additional regularity on the quantities involved and we do not have to specify the shape of the weight: it solves an inequality, different from \cite{CP82} where they impose the equality, and we just need to prove that the weight is bounded from above and below. Finally, thanks to Theorem \ref{the:deltaP_G0_Linfty}, we are able to pass to the incompressible limit for all fields of pressure.

\subsection{The Aronson-B\'enilan Estimate and Regularity Theory of the Porous Medium Equation}

In \cite{CF79}, the authors prove that the unique generalised solution to the porous medium equation in two (or more) dimension, a result due to \cite{Sab61}, is continuous. Moreover, they give an explicit expression for the modulus of continuity in space and time, \cf \cite[Theorem 1.1]{CF79},  thus extending the known one-dimensional result on the H\"older-regularity in space by Aronson, \cf \cite[Theorem, p.465]{Aro69}. It is important to stress that the regularity theory of Caffarelli and Friedman heavily relies on the AB-estimate. On the one hand, it allows them, in some sense, to quantify how the density at the centre of small balls  changes in small time instance, \cf \cite[Lemmas 2.2 \& 2.3]{CF79}, and derive an explicit modulus of continuity for the multi-dimensional porous medium equation. On the other hand, it allows them to study the regularity of the free boundary in any dimension, \cf \cite{CF80}, generalising the one-dimensional results of \cite{Aro70a, CF79a}. In doing so, they need to quantify how fast the density begins to intrude a previously unoccupied domain, \cf \cite[Lemma 2.1]{CF80}, which uses the multi-dimensional AB-estimate. One of the key findings in their paper is the H\"older-regularity of the free boundary, $t=f(x)$, which allows them to improve the modulus of continuity solutions which are established to be H\"older-continuous, too, \cite[Section 4]{CF79}. Eight years later, in \cite{CVW87}, the authors prove that the free boundary is in fact Lipschitz continuous for all times larger than the first time that the solution contains a ball which includes the initial data. For smaller time instance, however, H\"older regularity is optimal due to so-called focusing phenomena, \cite{AG93, AA95, CVW87}. The focusing problem is dedicated to understanding how areas of vacuum in the initial data are filled by the evolution of the porous medium equation. While solutions are H\"older continuous, the pressure gradient can blow up in time, \cf \cite{AG93}.

A hard problem is to perform a linear limit of the PME, to recover the heat equation. Indeed, the mentioned H\"older regularity for the solution of the PME cannot be obtained. Only in recent works \cite{G2017, GST2019}, the authors prove that given  an initial data with low regularity, the solution is $L^p$ in time and belongs to a fractional Sobolev space in $\R^d$ with $d>1$. Moreover, using the Barenblatt solution, the authors prove that this result is the optimal one.

\vspace{1cm}
\subsection{Applications to Modelling Tissue Growth}
In the last few decades, the study of cancer development improve due to new analytical tools and due  to the introduction of new numerical methods \cite{bellomo2003modeling, friedman2004hierarchy, bellomo2008foundations, lowengrub2009nonlinear, deroulers2009modeling, bellomo2009complexity}. The main difficulty in studying  these phenomena is the vast biological complexity related to the presence of different kinds of inter-specifically and intra-specifically interacting cells. Describing the tumour at the macroscopic level, we can distinguish two categories. On the one hand, we can devise partial differential equations (PDEs) of type , used primarily to model tumour growth \cite{greenspan1976growth, byrne1996modelling, ribba2006multiscale, preziosi2009multiphase, BCGRS10, ciarletta2011radial}, in which cells are represented by densities. On the other hand, tissue growth can be described by devising a free boundary model \cite{greenspan1972models,cui2008asymptotic,friedman2008stability}, where tissue growth is due to the motion of its boundary. Each of these approaches have their advantages: the first approach, also called mechanical models, is widely studied with many numerical and analytical tools. Regarding the second approach, it is closer to the biological vision of the tissue and allows to study its motion and dynamics. There is a well-developed technique to establish a link between the two approaches,  the so-called incompressible limit, which implies that the pressure becomes stiff \cite{PQV14,PV15, HV17,MPQ17,BPPS19,CDHV19,LTWZ19,KP15}.

\subsubsection{Population-based Description of Tissue Growth}
The simplest way to model tumour growth is introducing a single equation describing the evolution of the abnormal cell density, $n(x,t)$, where $x \in \R^d$ and $t \in \R^{+}$ which evolves under pressure forces and cell multiplication according to the equation \cite{BD09,PQV14}
\begin{align}
\label{eq:model}
    \frac{\partial n}{\partial t} -\diver(n \nabla p) = n G(p),
\end{align}
where $p=p(n)$ is the pressure field and $G:= G(p)$ models the proliferation of cells and it is called growth function. Suitable assumptions have to be imposed on $G$, for instance
\begin{align}
    \label{eq:AssumptionsOnG}
    G'(p) < 0, \qqand G(p_M) = 0,
\end{align}
which imply that the increasing number of the cells is limited by the pressure $p$ and $p_M >0$ is called the {\em homeostatic pressure} \cite{BD09}. Concerning the pressure, in many papers \cite{PQV14, PV15, HV17, BPPS19, LTWZ19,KP15}, there is an explicit and assigned relation between $p$ and $n$. In this work, to be as general as possible, we only assume that
\begin{align}
    \label{eq:pressure_function}
    p = p(n), \qquad p(0)=0, \qquad p'(n) >0,
\end{align}
for $n>0$. We have two examples in mind, the classical Power-law case, where 
\begin{align}
\label{eq:PowerLawPressure}
    p(n) = n^\gamma,
\end{align}
see \cite{PQTV14, BPPS19, david2020free} and the pressure used in \cite{DHV20} (called DHV throughout) where 
\begin{align}
    \label{eq:SingularPressure}
     p(n) = \varepsilon\frac{n}{1-n},
\end{align}
\cf \cite{DHV20}.
For a general  pressure law, the quantity $p$ satisfies the evolution equation given by
\begin{align}
\label{eq:p_growth}
    \partial_t p = \abs{\nabla p}^2 + q w,
\end{align}
where
\begin{align}
    \label{eq:DefOf_w}
  q(p) := np'(n), \qqand   w := \Delta p +G(p).
\end{align}
The aforementioned examples give 
\[
    q(p) = \gamma p   \quad \text{(polytropic law)}, \qqand q(p) =p\prt*{1+\frac{p}{\epsilon}}  \quad \text{(DHV law)}.
\]
They differ deeply near $p=0$ in their behaviours as $\epsilon \to 0$ and $\gamma \to \infty$ and this is a major issue if one wants to study these limits and establish the Hele-Shaw free boundary problem.

For $G \equiv 0$, Aronson and B\'enilan build their estimate on the observation that one can obtain an equation for $w$ \cite{AB79}. Their argument can be extended to include $G$ and leads to the equation
\[
    \partial_t w= 2 \sum_{i,j=1}^N (\partial^2_{ij} p)^2 + 2 \nabla p \cdot \nabla w - G'(p) |\nabla p|^2 + \Delta (qw) + G'(p) q w .
\]
Because
$$
    2 \sum_{i,j=1}^N (\partial^2_{ij} p)^2 \geq \frac{2}{N}(\Delta p)^2 = \frac2N (w-G)^2,
$$
and since we assume $G' \leq 0$, \ie, Eq. \eqref{eq:AssumptionsOnG}, we may also write 
\begin{align}
    \label{eq:Deltaw}
    \partial_t w \geq  \frac2N  w^2 + 2 \nabla p \cdot \nabla w  + \Delta (qw) + \brk*{G'(p)q- \frac4N G } w .
\end{align}
This inequality, which is self-contained when $G\equiv 0$, is the very basis of estimates on $|w|_- = \max(0, -w)$ which we analyse in different $L^p$ spaces.
\\

\subsubsection{Free Boundary-based Description of Tissue Growth}
Besides its huge impact on the regularity theory of solutions to the porous medium equation and the free boundaries thereof, the AB-estimate proves to be a crucial tool for building a bridge between a density-based description and a geometric description of tissue growth. The link between the two models is established through a rigorous  study of the \emph{incompressible limit} of the porous medium pressure equation, \cf Eq. \eqref{eq:p_growth}, as the pressure law becomes stiffer and stiffer, \ie, $\gamma \to \infty$ or  $\epsilon\to0$,  in the respective pressure law (Eqs. \eqref{eq:PowerLawPressure}, \eqref{eq:SingularPressure}). As a result an incompressible model is obtained, satisfying  two  relations. The first,  $p (n-1) = 0$, implies the absence of any pressure in zones that are not saturated ($\{n<1\}$), while the second one, also referred to as \emph{complementarity relation}, yields an equation satisfied by the pressure on $\{p>0\}$, which is of the form
$$
    p (\Delta p + G(p)) =0.
$$ 
It is immediately apparent that strong regularity  is needed to obtain such an expression, which is provided by (adaptations) of the AB-estimate --- bounds on the Laplacian of the pressure are enough to infer strong compactness of the pressure gradient. This was first observed in \cite{PQV14} to be equivalent to being able to pass to the limit in the porous medium pressure equation and obtain the incompressible limit.

\subsubsection{Extension to Two Species}
Let us highlight that the mathematical theory of the limit for  equations like Eq. \eqref{eq:model} is well studied both with $G\neq 0$, \cf \cite{PQV14, DHV20, KP15}, and without $G$, \cf  \cite{benilan1989limit, gil2001convergence, gil2003boundary}, as well as in the case where nutrients and viscosity are included, \cf  \cite{PV15,  DS20, DPSV20, david2020free}. The limit model turns out to be a free boundary model of Hele-Shaw type.

The model of a single evolution equation, which describes the tumor cell distribution, can be complemented by another species consisting of healthy tissue, and it is given by
\begin{align}
    \label{eq:two_species}
    \frac{\partial n^{(i)}}{\partial t} - \diver\prt*{n^{(i)} \nabla p} = n^{(1)} F^{(i)}(p)+  n^{(2)} G^{(i)}(p),
\end{align}
where $i = 1,2$, $n^{(1)}, n^{(2)}$ denote the population densities and $G^{(i)}, F^{(i)}$ model the reaction or growth phenomena, which are assumed to depend exclusively on the pressure according to experimental observations \cite{BD09,RBEJPJ10}. The  system structure of Eq. \eqref{eq:two_species} causes serious analytical difficulties, \cf \cite{CFSS17,GPS19,BPPS19, PX20, DJ20}, due to its hyperbolic flavour. The careful study of the pressure equation helps in proving the existence of solutions and obtaining uniform estimates with respect to the stiffness parameter, \ie, $\gamma$ in the classical Power-law case, $p = n^{\gamma}$ \cite{PQTV14,MPQ17,GPS19}. Due to the insufficient regularity of the pressure, the incompressible limit can be achieved just in $1D$ using Sobolev embedding. For the DHV pressure law, \ie,  assuming that the pressure blows up at a finite threshold, \cf \cite{HV17,CDHV19,DHV20}, similar mathematical difficulties arise: in order to pass to the limit, strong restrictions have to be imposed.

While the incompressible limit for multiple species remains an interesting open problem for the Darcy law, including viscosity of cells in the model, \ie, altering the velocity, $u$, in Eq. \eqref{eq:continuity_eqn}, changes the analytical properties of the model drastically and recently. In case of the so-called {\em Brinkman’s law} \cite{A99_brinkman}, the model reads
\begin{equation}
    \label{eq:two_species_brinkman}
    \left\{
    \begin{aligned}
    &\frac{\partial n^{(i)}}{\partial t} - \diver(n^{(i)}\nabla W) = n^{(1)} F^{(i)}(p)+  n^{(2)} G^{(i)}(p)\\
    & - \nu \Delta W + W = p,
    \end{aligned}
    \right.
\end{equation}
where $i = 1,2$ and $\nu$ is the viscosity parameter. On the one hand, the idea to couple the two equations for the individual species through Brinkman’s law changes the behaviour dramatically as mentioned before, \ie, classical techniques used in the one-specie model fail \cite{PV15}, but on the other hand the pressure field gains regularity and some mathematical difficulties can be overcome. Recently, two results have been obtained in this direction: in \cite{DS20}, the authors are able to establish the incompressible limit in the one dimensional case by establishing uniform BV-bounds for the two species; then in \cite{DPSV20}, since the BV-strategy fails in higher dimensions, by employing a non-local compactness criterion \cite{J10}, the passage to the incompressible limit can be accomplished.

\bigskip

\subsection{Plan of the paper}
The paper is organised as follows. In Section \ref{sec:L1}, we derive the $L^1$-estimate  without and with a reaction term $G$ and we show the advantages and limitations of adding a weight $h$. Then, in Section \ref{sec:Linfty}, we conduct similar computations calculations in the $L^\infty$ space proving that the additional weight helps in generalising the AB-type estimate for all pressure laws. In Section \ref{sec:L2}, we perform an $L^2$-estimate showing that it can be closed just for the a particular class of pressure laws. Finally, in Section \ref{sec:conclusions}, we add few concluding remarks.
\newpage

\section{$L^1$-type estimate}
\label{sec:L1}
The form of Eq. \eqref{eq:Deltaw} is well adapted to perform $L^1$-estimates of the second-order quantity $w$ because it generates, thanks to the Kato inequality, the following inequality
\begin{align} \label{eq:Deltaw1}
\partial_t |w|_- \leq - \frac2N  |w|_-^2 + 2 \nabla p\cdot \nabla |w|_- 
 + \Delta (q |w|_-) + \left[G'(p)q- \frac4N G \right] |w|_- .
\end{align}
It turns out that further manipulations lead to restrictions which are more demanding than expected. To explain that, we first treat  the case $G\equiv 0$. The extension to the case $G \geq 0$ is stated in the second corresponding subsection.

\subsection{$L^1$-Estimates when $G\equiv 0$}
$\mbox{}$\\

When $G \equiv 0$, $w = \Delta p$ and a simple integration of Eq. \eqref{eq:Deltaw} yields,  
\[
\ddt  \int_{\R^N} | w |_- \dx{x} \leq 2 \frac {N-1}{N} \int_{\R^N} | w |_-^2\dx{x}.
\]
Because of the quadratic growth of the right-hand side,  this inequality provides us with an $L^1$ control only in dimension $N=1$.  We thus adopt a different strategy. By adding a positive weight function, $h= h(p)$, we aim to study whether or not it  helps improve the above result. We shall establish the

\begin{theorem}[Case $G\equiv 0$, $L^1$-theory] Assume the pressure law is such that for $p>0$
\begin{equation}
\label{AS-L1}
    \alpha_1(p) := \int_0^p \left[ q(\cdot)h'(\cdot) + \frac {1}{N}h(\cdot)\right] \dx{\rho}>0, \quad \text{with } \quad 
    h(p):=\int_0^p  e^{\int_\cdot^\rho \frac{1}{q}} \; \dx{\rho}.
\end{equation}   
Then, the following a priori estimates hold true 
\begin{equation}
    \label{estL1-1}
    \int_{\R^N} h(p) |\Delta p(t)|_- \dx{x} + 2\int_0^t \int_{\R^N} \alpha_1 (p(s)) |\Delta p (s)|_-^2 \dx{x}\dx{s} \leq \int_{\R^N} h(p) |\Delta p(0)|_-\dx{x} ,
\end{equation}
and for all $t \leq T$
\begin{equation}
    \label{estL1-2}
    \int_{\R^N} h (p) |\Delta p(t)|_- \dx{x} \leq \frac {A(T)}{t} \qquad \text{with } \qquad  A= \frac 12 \sup_{0\leq t \leq T} \int_{\R^N} \frac{h(p)^2}{\alpha_1(p)} \dx{x}.
\end{equation}
Estimate \eqref{estL1-1} also holds when $N=1$ with $h\equiv 1$ and $\alpha_1\equiv 0$.
\label{th:L1G=0}\end{theorem}

Notice that, because $q'(p) p'(n)= p'(n)+n p''(n)$, the condition $\alpha_1 >0$ is satisfied when, for instance,  $p(n)$ is a convex function.

However, it turns out that the result of this theorem is rather weak compared of $L^\infty$-type estimate. Indeed, the calculation forces us to choose $h(0)=0$, therefore  the estimate is weak near the free boundary $\{p=0\}$. However, for the pressure laws we have introduced before, when $\gamma$ is large or $\epsilon$ is small, the decay rate scales correctly with $\gamma$ or $\epsilon$, and allows for a uniform control of $\partial_t p$. To see this we may give the corresponding expressions of $\alpha_1, \; h, \; A$, explicitly:

For the Power-law, we have $q(p)=\gamma p$ and 
\[
h(p)= \frac{\gamma}{\gamma+1} p^{\frac{\gamma+1}{\gamma}}, \qquad
\alpha_1 (p) = \gamma \left[1-  \frac {N-1}{N}  \frac {1}{\gamma+1}\right]p^{\frac{\gamma+1}{\gamma}}>0,
\qqand A =\mathcal{O}(\gamma^{-1}) .
\]
\\
For the DHV law, the situation is similar and we find 
\[
h(p)=\epsilon p - \epsilon^2 \ln\left(1+ \frac p \epsilon\right), \quad
\alpha_1 (p) = p^2 + \epsilon \frac {N-1}{N} \left[ \epsilon \ln\left(1+ \frac p \epsilon\right) -p \right]>0,\qand A =\mathcal{O}(\epsilon^2).
\]

\begin{proof}
Hence, integrating over $\R^N$, we get
\begin{align}
\label{eq:lastvalue_L1h}
\begin{split}
    \ddt \int_{\R^N} h(p) \abs{w}_{-} \dx{x} 
    &= \int_{\R^N}  \prt*{h'(p)\, \abs{w}_{-} \prt*{\frac{\partial p}{\partial t}} + h(p)\, \prt*{\frac{\partial \abs{w}_{-}}{\partial t}}}\dx{x}
    \\
    &\leq  \int_{\R^N} h'(p)\, \abs{\nabla p}^2 \abs{w}_{-} \dx{x} -\int_{\R^N} h'(p)\, q(p)\, \abs{w}_{-}^2 \dx{x} - \int_{\R^N} \frac{2 h}{N} \abs{w}_{-}^2 \dx{x}  
    \\
    &\quad + \underbrace{\int_{\R^N} h(p)\, \Delta \prt*{q(p)\, \abs{w}_{-}} \dx{x}}_{I_1} + \underbrace{2\int_{\R^N} h(p)\, \nabla p \cdot \nabla \abs{w}_{-}\dx{x}}_{I_2}.
\end{split}
\end{align}
We next estimate separately the integral terms $I_1$ and $I_2$. Beginning with $I_1$, integrating by parts twice, we obtain
\begin{align*}
    I_1  &= \int_{\R^N} q(p)\, \abs{w}_{-} \Delta h(p) \dx{x} 
    \\
    &= \int_{\R^N} q(p)\, h''(p) \, \abs{\nabla p}^2 \abs{w}_{-} \dx{x} - \int_{\R^N} h'(p)\, q(p)\, \abs{w}_{-}^2 \dx{x}.
\end{align*}
As for the second term, $I_2$, we integrate by parts once, which yields
\begin{align*}
    I_2 &= - 2 \int_{\R^N} h'(p) \abs{\nabla p}^2 \abs{w}_{-}\dx{x} - 2 \int_{\R^N} h(p)\, \abs{w}_{-}^2 \dx{x}.
\end{align*}
Substituting the expressions for the two integral terms $I_1$ and $I_2$ into Eq. \eqref{eq:lastvalue_L1h}, we obtain
\begin{align}\label{eq:L1G=0}
    \ddt \int_{\R^N} h(p)\, \abs{w}_{-} \dx{x}
    \leq - 2\int_{\R^N} \alpha_1 \, \abs{w}_{-}^2\dx{x} + \int_{\R^N} \beta_1\,  \abs{w}_{-} \abs{\nabla p}^2 \dx{x},
\end{align}
where
\begin{align}
    \label{eq:ConditionOnAlpha1Beta1}
    \alpha_1 := {h'q}-  h \prt*{1-\frac{1}{N}}, \qqand \beta_1  :=q h'' - h'.
\end{align}

In dimension $N=1$ we can choose $h=1$, $\alpha_1 = \beta_1=0$ which provides us with and $L^1$-estimate of $|\Delta p|_-$ and proves the last statement of Theorem~\ref{th:L1G=0}.

In higher dimension, we are unable to do that and we solve $\beta_1 (p)=0$ instead, which gives the expression of $h$ and $\alpha_1$ in \eqref{AS-L1}.
Then, integrating Eq.~\eqref{eq:L1G=0} in $t$ gives the announced estimate~\eqref{estL1-1}.
\\

To obtain estimate~\eqref{estL1-2}, we use the Cauchy-Schwarz inequality to write
\[
    \left(\int_{\R^N} h \abs{w}_{-}\right)^2 \dx{x} \leq A \; \int_{\R^N} \alpha_1 \, \abs{w}_{-}^2\dx{x}, \qquad A:= \frac 12\int_{\R^N} \frac{h^2}{\alpha_1} \dx{x}.
\]
Substituting this information into Eq. \eqref{eq:L1G=0}, we get
$$
    \ddt \int_{\R^N} h(p)\, \abs{w}_{-} \dx{x}\leq - A^{-1} \left(\int_{\R^N} h \abs{w}_{-}\dx{x} \right)^2 .
$$
Since $U(t)= \frac A t $ is a solution, we conclude that
$$
    \int_{\R^N} h(p)\, \abs{w}_{-} \dx{x} \leq A /t,
$$
which proves the statement \eqref{estL1-2} and concludes the proof of Theorem~\ref{th:L1G=0}.
\end{proof}

\subsection{$L^1$-estimates with $G\neq 0$}
\label{sec:L1G}
$\hbox{}$\\

With the notations of Theorem \ref{th:L1G=0}, we define
\[
\bar \delta_1 = \max_{0\leq p \leq p_M}\left\{ G \left[ 2\frac{N-2}{N} -\frac{h'q}{h}\right]+ G'q \right\}.
\]
\begin{theorem}[Case $G' \leq 0$, $L^1$-theory] With the notations and assumptions of Theorem \ref{th:L1G=0}, $G' \leq 0$ and $\bar \delta_1$ as above, the following a priori estimates hold true 
\begin{equation}
    \label{estL1-G1}
    \int_{\R^N} h(p) |w(t)|_- \dx{x} + \int_0^te^{\bar \delta_1 (t-s)} \int_{\R^N} \alpha_1 (p) |w|_-^2 \dx{x}\dx{s} \leq e^{\bar \delta_1 t}\int_{\R^N} h(p) |w(0)|_-\dx{x}, 
\end{equation}
and
\begin{equation}
    \label{estL1-G2}
    \int_{\R^N} h (p) |w(t)|_- \dx{x} \leq A \frac{\bar \delta_1 e^{\bar \delta_1 t} }{e^{\bar \delta_1 t}-1} .
\end{equation}
with $A$ as in Eq. \eqref{estL1-2}.
Estimate \eqref{estL1-G1} also holds when $N=1$ with $h\equiv 1$ and $\alpha_1\equiv 0$.
\label{th:L1G}\end{theorem}

Notice that the sign of $\bar \delta_1$ does not play a role here.
\\

\begin{proof}

Still building on the inequality \eqref{eq:Deltaw} and  using a positive weight, $h= h(p)$, the evolution in time of the quantity $h w$ is given by
\begin{align*}
	\partial_t\left(h\abs{w}_{-}\right) &= h'\abs{w}_{-} \partial_t p+h\partial_t \abs{w}_{-} \\
	& \leq  h' \abs{w}_{-} \left(\abs{\nabla p}^2 +qw\right)  
	\\ &\qquad + h  \set*{
	- \frac2N  |w|_-^2+ 2 \nabla p\cdot \nabla |w|_- + \Delta (q |w|_-) + \left[G'(p)q- \frac4N G \right] |w|_- }
\end{align*}	
Integrating, we get
\begin{align*}
\begin{split}
    \ddt \int_{\R^N} h\abs{w}_{-}  \dx{x} &\leq \int_{\R^N} h'\abs{\nabla p}^{2}\abs{w}_{-}\, \dx{x} -  \int_{\R^N} h' q \abs{w}_{-}^2 \, \dx{x} -\int_{\R^N} \frac{2h}{N} \wm^2 \, \dx{x} 
    \\
    &\quad + \int_{\R^N} h \Delta \prt*{q\abs{w}_{-}}\, \dx{x} + 2 \int_{\R^N} h \nabla p \cdot \nabla\abs{w}_{-} \, \dx{x}  -\int_{\R^N} \left[ \frac{4 }{N} G- G' q\right] h \; \wm  \, \dx{x} 
    \\[1em]
    &\leq \int_{\R^N} h'\abs{\nabla p}^{2}\abs{w}_{-}\, \dx{x} -\int_{\R^N}  \left[\frac{2h}{N}+  h' q \right]\, \wm^2 \dx{x} - \int_{\R^N} \left[ \frac{4 }{N} G- G' q\right] h \; \wm  \, \dx{x}   \\
    &\quad + \underbrace{\int_{\R^N}  q\abs{w}_{-} \; \Delta h \, \dx{x}}_{I_1} -2\;  \underbrace{ \int_{\R^N} \left[\nabla h\cdot \nabla p + h \Delta p\right]\abs{w}_{-} \, \dx{x}}_{I_2} .
    \end{split}
\end{align*}
Integrating by parts twice, the term $I_1$ can be rewritten as
\begin{align*}
   I_1 &= \int_{\R^N} \left[h'' \abs{\nabla p}^2 +h'\Delta p\right] q \abs{w}_{-} \, \dx{x}\\
    &= \int_{\R^N} h'' q \abs{\nabla p}^2 \abs{w}_{-} \, \dx{x} -\int_{\R^N} h' q\wm^2\, \dx{x} -  \int_{\R^N} h'qG\abs{w}_{- }\, \dx{x}.
\end{align*}
Next we simplify the term $I_2$. Integrating by parts and using the chain rule we obtain
\begin{align*}
	I_2 &=  \int_{\R^N} h'\abs{\nabla p}^2 \abs{w}_{-} \, \dx{x} -  \int_{\R^N} h\abs{w}_{-}^2 \dx{x} -   \int_{\R^N} h G\abs{w}_{-}\, \dx{x}.
\end{align*}
Substituting $I_1$ and $I_2$ back into our main inequality, we get
\begin{align*}
    \ddt \int_{\R^N} h \abs{w}_{-} \dx{x}
    &\leq - 2\int_{\R^N} \alpha_1^h\abs{w}_{-}^2 \, \dx{x} + \int_{\R^N} \beta_1^h \abs{w}_{-} \abs{\nabla p}^2\, \dx{x} +\int_{\R^N} \delta_1^h h \abs{w}_{-} \, \dx{x},
\end{align*}
where
\begin{align*}
    \alpha_1^h = {h'q} - h\frac{N-1}{N}, \qquad \beta_1^h = {h'' q}-{h'}, \qquad\text{and}\qquad \delta_1^h = G\left[ 2\frac{N-2}{N} -\frac{h'q}{h}\right]+ G'q \leq \bar \delta_1 .
\end{align*}
To control the terms on the right-hand side, we argue as in the case $G\equiv 0$. In dimension $N=1$, we can choose $h=1$, $\alpha_1^h=\beta_1^h=0$. Otherwise, it is sufficient to impose that
\[
    \alpha_1^h (p) >0 , \qqand \beta_1^h (p)=0,
\]
as chosen in Theorem \ref{th:L1G=0}.

We rewrite the inequality, after using the Cauchy-Schwarz inequality, as 
\begin{align*}
    \ddt  \int_{\R^N} h \abs{w}_{-} \dx{x}
    &\leq - A^{-1} \left(\int_{\R^N} h\abs{w}_{-}\right)^2 \, \dx{x} + \bar \delta_1 \int_{\R^N} h \abs{w}_{-} \, \dx{x}.
\end{align*}
It remains to observe that $A \frac{\bar \delta_1 e^{\bar \delta_1 t} }{e^{\bar \delta_1 t}-1}$ is a solution and we obtain the statements of Theorem \ref{th:L1G}.
\end{proof}

\bigskip

\section{$L^\infty$-type Estimate}
\label{sec:Linfty}
The other extreme $L^p$-space is that used in the original paper, and establishes a bound in $L^\infty$ of $|w|_-$. It uses the strong form of the equation satisfied by $w$. Namely, starting from Eq. \eqref{eq:Deltaw}, we can write
\[
    \partial_t w \geq  \frac2N  w^2 + 2 \nabla (p+q) \cdot \nabla w + q \Delta w + w \Delta q + \brk*{G'(p)q - \frac4N G } w .
\]
Using that 
\[
    \Delta q =q''(p) |\nabla p|^2 + q'(p) \Delta p =  q''(p) |\nabla p|^2 + q'(p) (w-G),
\]
we find 
\begin{align}
    \label{eq:DeltawStrong}
    \partial_t w \geq  \left[q'+\frac2N \right]  w^2 + 2 \nabla (p+q) \cdot \nabla w +q'' |\nabla p|^2 w + q \Delta w + \brk*{G'\; q - \left(\frac4N+q'\right) G } w .
\end{align}

A first result that can be deduced directly from this calculation is
\begin{theorem}[Lower bound on $\Delta p$, special case]
    \label{the:max_principle_special}
    Assume that $q''\geq 0$ and assume there are constants  $\bar \delta \in \R$, $\alpha_0>0$ such that  $G'\; q - (\frac4N+q') G \leq \bar \delta$ and $q'+\frac2N \geq \alpha_0$, then we have
\[
     \Delta p + G\geq -\frac{1}{\alpha_0} \frac{\bar \delta e^{\bar \delta t} }{e^{\bar \delta t}-1}.
\]
\end{theorem}

This result applies to homogeneous  pressure laws $q(p) = \gamma p$ and $\alpha_0=O(\gamma)$
and matches that of \cite{PQV14}, see \cite[Eqs. (2.14)]{PQV14}. It also applies to DHV law, $q = p+\varepsilon^{-1}p^2, q' = 1 + 2\varepsilon^{-1}p$ and $q'' = 2\varepsilon^{-1}$, but then $\alpha_0=O(1)$ (when $p\approx 0$) does not give a uniform decay as $\frac 1 \epsilon)$ as needed to study the Hele-Shaw limit.
\\

To treat more general pressure laws we can refine the argument and as in the $L^1$-case, we begin with the porous-medium equation and then we include the growth term.

\subsection{$L^\infty$-estimates when $G\equiv 0$}
$\mbox{}$\\

We begin with estimating the Laplacian of the pressure, $w = \Delta p$, when $G\equiv 0$. We are going to prove the following theorem
\begin{theorem}[Lower bound on $\Delta p$, $G\equiv 0$] 
\label{the:deltaP_G0_Linfty}
1. Assume that
\[
\left(\frac q p\right)' \geq 0, \qqand \widetilde \alpha_0:= \min_{0\leq p \leq p_M} \frac{p \; q'(p)}{q(p)} >0, 
\]
then
\[
\frac {q(p)}{p} \;  \Delta p \geq - \frac {1}{\widetilde \alpha_0 \; t}.
\]
2. Assume that $q'(p) >-1$ for $p\in [0,p_M]$, and
\[
\widetilde \alpha_0=  \frac{\min _{0 \leq p \leq p_M} q'(p)} {1+ \max_{0 \leq p \leq p_M} q'(p)}+ \frac 2 N \frac{1+q'(0)} {\left(1+ \max_{0 \leq p \leq p_M} q'(p)\right)^2} >0.
\]
Then, we have
\[
\left(1+ \min_{0 \leq p \leq p_M} q'(r)\right)   \Delta p \geq - \frac {1}{\widetilde \alpha_0 \; t}.
\]
\label{th:Linfty_G=0}
\end{theorem}
 With the first set of assumptions, the estimate is compatible with the Hele-Shaw asymptotics in the two examples of Power-Law (then $\widetilde \alpha_0= 1$ and $\frac q p = \gamma$) and DHV law (then $\widetilde \alpha_0= 1$ and $\frac q p = \frac {p+\epsilon}{\epsilon}$). The second set of assumption is an explicit example motivated by \cite{CP82}.
\\

To understand if there is some slack in the estimate, we compute the evolution of the quantity $h w$, where $h = h(p)$ is assumed to be a positive weight function. We obtain
\begin{align*}
        \frac{\partial (h w)}{\partial t} &= w h' \frac{\partial p}{\partial t} +h \frac{\partial w}{\partial t}\\
        &\geq  w h' \prt*{\abs{\nabla p}^2+qw}
        \\
        &\quad + h \left(\left[q'+\frac2N \right]  w^2 +q'' |\nabla p|^2 w + 2 \nabla (p+q) \cdot \nabla w + q \Delta w  \right),
\end{align*}
where we used the equation satisfied by the pressure, Eq. \eqref{eq:p_growth},  and Eq. \eqref{eq:DeltawStrong} for $w$. Upon rearranging the terms, we get
\begin{align}
    \label{eq:LinftyFor_hv}
     \frac{\partial (h w)}{\partial t}
       & \geq w^2 \brk*{{hq'}+ qh'+\frac2N} + w \abs{\nabla p}^2\prt*{{h'}+ q''h}+ q\underbrace{h \Delta w}_{I_1}+\underbrace{2h (1 + q') \nabla p \cdot \nabla w}_{I_2} .
\end{align}
The terms inwolving the linear operators have to be rewritten in terms of the new quantity $hw$ rather than $w$. Therefore, from the first term we get 
\begin{align*}
    I_1 &=  \Delta\prt*{hw}-w \Delta h - 2\nabla h\cdot \nabla w \\
    & =  \Delta \prt*{hw} - w \prt*{h''\abs{\nabla p}^2 + h' w} - 2h' \nabla p \cdot \nabla w\\
    &= \Delta \prt*{hw} - w \prt*{h''\abs{\nabla p}^2 + h' w} - 2\frac{h'}{h} \nabla p \cdot \prt*{\nabla\prt*{hw}- w \nabla h}\\
    &=  \Delta \prt*{hw} -2\frac{h'}{h} \nabla p \cdot \nabla\prt*{hw}  - w \prt*{h''\abs{\nabla p}^2 + h' w} + 2 \prt{hw} \frac{\prt{h'}^2}{h^2} \abs{\nabla p}^2\\
    &=  \Delta \prt*{hw} -2\frac{h'}{h} \nabla p \cdot \nabla\prt*{hw} + w \abs{\nabla p}^2 \prt*{2\frac{\prt*{h'}^2}{h}-h''} - h' w^2.
\end{align*}
For the second one we have
\begin{align*}
    I_2 &= 2  \nabla (p+q) \cdot h\nabla w \\
    &=  2 \nabla (p+q) \cdot \prt*{\nabla \prt*{hw}-h'w\nabla p}\\
    &= 2 \nabla (p+q) \cdot \nabla \prt*{hw} - 2 \prt*{1+q'}h' w \abs{\nabla p}^2 .
\end{align*}
Substituting the simplified expressions of $I_1$ and $I_2$ into Eq. \eqref{eq:LinftyFor_hv}, we have
\begin{align} \label{eq:deltapLinfinity}
    \frac{\partial (h w)}{\partial t}  \geq \alpha_{\infty}^h\prt*{hw}^2 + \beta_{\infty}^h (hw) \abs{\nabla p}^2 + \curlyL^h_{\infty}(hw),
\end{align}
where 
\begin{align*}
    \alpha_{\infty}^h = \frac{q'}{h}+\frac{2}{Nh^2} , \qqand \beta_{\infty}^h =  q''-\frac{qh''}{h} +2\frac{q\prt*{h'}^2}{h^2}-2\frac{q'h'}{h}-\frac{h'}{h},
\end{align*}
as well as
\begin{align*}
    \curlyL^h_{\infty}(hw) = q \Delta \prt*{hw}+  2 \prt*{1+q'-\frac{qh'}{h}} \nabla p \cdot \nabla\prt*{hw}.
\end{align*}

In order to find a sub-solution and to close the estimate, it is enough to ensure that
\begin{equation}
    \label{cond1}
     \alpha_{\infty}^h = \frac{q'}{h}+\frac{2}{Nh^2}  \geq \al_0 >0, \qquad  \frac{ \beta_{\infty}^h}{h} = \left(\f{q}{h}\right)'' + \left(\f{1}{h}\right)'\leq  0,
\end{equation}
where $\alpha_0$ is a constant. We propose two strategies to fulfill these requirements.

\paragraph{1. Assume $\left(\frac q p\right)' \geq 0$.}
With this assumption, we can simply choose 
\[
h=\frac{q}{q'(0)p}, \qquad \left(\f{q}{h}\right)' = q'(0), 
\]
because $\frac 1h$ is non-increasing. Then we compute, 
\[
\alpha_0:=q'(0)\simeq \min_{0\leq p \leq p_M} \left[ \frac{p \; q'(p)}{q(p)}+ \frac 2N \frac{q'(0)p^2}{q(p)^2}\right] ,
\]
and this gives our first statement in Theorem~\ref{th:Linfty_G=0} after simplifying a coefficient $q'(0)$.

\paragraph{2. Assume $q'(p) > -1$.}
Then we impose, from Eq. \eqref{cond1}
\begin{align}
    \label{eq:diff_equation}
 \left(\f{q}{h}\right)' + \f{1}{h}=  1+q'(0), \qqand h(0)=1.
\end{align}
Because of the degeneracy at $p=0$, the condition $h(0)=1$ is imposed since we can compute
\[
    1+q'(0)=\frac{q'(0)}{h(0)} + q(0) (h^{-1})'(0) + \frac{1}{h(0)} = \frac{q'(0)+1}{h(0)}.
\]
We can analyse the differential equation \eqref{eq:diff_equation}.

\begin{lemma}
\label{lemma:cond_Linfty}
Assume that $q'(p) > -1$ for $p\in [0, p_M]$. 
Then, the solution of \eqref{eq:diff_equation} satisfies, 
\begin{align*}
   \frac{1+ \min_{0 \leq r \leq p} q'(r)} {1+q'(0)} \leq h(p) \leq  \frac {1+ \max_{0 \leq r \leq p} q'(r)}{1+q'(0)}, \qquad \forall p\in (0, p_M). 
\end{align*}
\end{lemma}
\begin{remark}[Power-law, DHV Pressure, General Pressure Laws]
	Notice that if $q$ is smooth as in DHV, $q'(0)=1$ but for the Power-law, $q'(0)=\gamma$.
	As a matter of fact, for the Power-law pressure and the singular pressure we are able to provide explicit expressions for $h$, \ie,
    \begin{align}
         h_{\gamma}(p) = 1,\qqand 
         h_{\e}(p)  =  1+ p - \e  \log (p+\e )+ \e \log \e.
     \end{align}
    However, we also emphasise that the bounds established in the preceding theorem allow to prove an Aronson-Bénilan type estimate for more general pressure-laws $p = p(n)$, under rather week assumptions on $q$, thus extending the known cases.
\end{remark}

\begin{proof}
First, we change variables and set, in Eq. \eqref{eq:diff_equation},
\begin{align*}
    u(p):= \frac1{h(p)},
\end{align*}
In the new variable, it becomes 
\begin{align}
\label{eq:qu_condition}
    qu' +(q'+1) u = 1+ q'(0),
\end{align}
with $ u(0)=1$. The rest of the argument, \ie, the proof of the upper and lower bound on $h$, is by contradiction. To this end, we define 
\[
    U(p) = a \; \frac{1+q'(0)} {1+ \max_{0 \leq r \leq p} q'(r)}.
\]
By construction, it is a non-increasing function and it satsifies $u(0)=1> U(0)=a$. Assume there exists a point $p^* \in (0,p_M)$, that we can choose to be minimal, such that
\begin{align}
    \label{eq:u_below_p_star}
    u(p^*)  = U(p^*).
\end{align}
Therefore, at this point we have to have  $u'(p^*)\leq 0$, as well as  $\big(1+q'(p^*)\big) u(p^*)=\big(1+q'(p^*)\big) \; U(p^*)$. Revisiting Eq. \eqref{eq:qu_condition}, we see that
\begin{align*}
      1+q'(0)  &=   (qu)'(p^*)+u(p^*)\\
      &=  q(p^*)u'(p^*) + (q(p^*) + 1)u(p^*)\\
      &\leq (1 + q(p^*))U(p^*)\\
      &\leq a (1+q'(0)),
\end{align*}
having used the fact that $qu'(p^*) \leq 0$ and the definition of $U$. It is clear that this is a contradiction, since $0<a<1$ and proves that $u > a U$. Finally, taking $a\to 1$, we obtain the upper bound. The lower bound is obtained in the same way.
\end{proof}

The derivation of Theorem~\ref{th:Linfty_G=0} is now as usual because $-\frac{1}{\alpha_0 t}$ is a sub-solution of Eq.~\eqref{eq:deltapLinfinity} and, using Lemma~\ref{lemma:cond_Linfty}, we can choose
$$
\alpha_0= \min _{0 \leq p \leq p_M} q'(p) \frac{1+q'(0)} {1+ \max_{0 \leq p \leq p_M} q'(p)}+ \frac 2 N \left(\frac{1+q'(0)} {1+ \max_{0 \leq p \leq p_M} q'(p)}\right)^2.
$$
Then, using Lemma~\ref{lemma:cond_Linfty} a second time,
\[
    \frac{1+ \min_{0 \leq r \leq p} q'(r)} {1+q'(0)} \;w \geq\min_p h(p) \;  w \geq -\frac{1}{\alpha_0 t} ,
\]
which gives the result of Theorem~\ref{th:Linfty_G=0}.

\subsection{$L^\infty$-estimates when $G\neq 0$}
$\hbox{}$\\

Next, we proceed by incorporating reaction terms, \cf Eq. \eqref{eq:model} and prove the

\begin{theorem}[Lower bound on $\Delta p$, general $G$] With the assumptions and notations of Theorem~\ref{the:deltaP_G0_Linfty},   we have, with the constant $\delta_\infty^h$ defined below
\[
    \left(1+ \min_{0 \leq p \leq p_M} q'(r)\right)   \left(\Delta p  + G \right)\geq - \frac {1}{\widetilde \alpha_0} \frac{\bar \delta_\infty^h e^{\bar \delta_\infty^h t} }{e^{\bar \delta_\infty^h t}-1}.
\]
\label{th:Linfty_G}
\end{theorem}

Still using \eqref{eq:DeltawStrong}, we compute the evolution of the quantity $hw$, \ie,
\begin{align}
\label{eq:ddt_hw}
\begin{split}
    \frac{\partial (h w)}{\partial t} 
    &\geq  h'\abs{\nabla p}^2 w + \prt*{\frac{h'q}{h^2}} \prt*{hw}^2 + hq'' \abs{\nabla p}^2 w+\prt*{\frac{q'}{h}}\prt*{hw}^2 \\[0.3em]
    &\quad +\underbrace{2  h \nabla (p+q) \cdot \nabla w}_{I_1} + q \underbrace{h \Delta w}_{I_2} + \prt*{hw} \prt*{ G'q- \prt*{\frac 4N +q'}G}.
\end{split}
\end{align}
Next, we rewrite the terms $I_1$ and $I_2$ using the new variable, $hw$. The first term, $I_1$, becomes
\begin{align*}
    I_1 &=  2 \nabla (p+q) \cdot \nabla \prt*{hw} - 2 w \abs{\nabla p}^2 \prt*{1+q'}h'.
\end{align*}
Regarding the second one, we obtain
\begin{align*}
    I_2 
    &=  \Delta\prt*{hw}-w \Delta h - 2\nabla h\cdot \nabla w \\
    & =  \Delta \prt*{hw} - w \prt*{h''\abs{\nabla p}^2 + h' \Delta p} - 2h' \nabla p \cdot \nabla w\\
    &= \Delta \prt*{hw} - w \prt*{h''\abs{\nabla p}^2 + h'w - h'G} - 2\frac{h'}{h} \nabla p \cdot \prt*{\nabla\prt*{hw}- w h' \nabla p}\\
    &= \Delta \prt*{hw} -2\frac{h'}{h} \nabla p \cdot \nabla\prt*{hw}  + w \abs{\nabla p}^2 \prt*{2\frac{h'^2}{h}- h''} -\prt*{hw}^2 \prt*{\frac{h'}{h^2}} + \prt*{hw}\prt*{\frac{h'G}{h}}.
\end{align*}
Substituting the terms $I_1$ and $I_2$ back into Eq. \eqref{eq:ddt_hw}, we get
\begin{align}
   \frac{\partial (h w)}{\partial t} &\geq   \alpha_\infty^h \prt*{hw}^2+\beta_\infty^h\abs{\nabla p}^2 w + \curlyL_\infty^h(hw) + \delta_\infty^h \prt*{hw} ,
\end{align}
where $\alpha_\infty^h$, $\beta_\infty^h$ and $\curlyL_\infty^h(hw)$ are  as  in inequality \eqref{eq:deltapLinfinity}, and 
\begin{equation}
	\label{cond3_GP}
    \delta_\infty^h(p) := G'q- \left(\frac 4 N +q'\right)G +\frac{qh'G}{h} \leq \bar \delta_\infty < \infty.
\end{equation}

In order to close the $L^{\infty}$-type bound, it suffices to apply a slight variation of the proof of Theorem  \ref{th:Linfty_G=0} incorporating additional terms related to the growth. Again, we have to require that
\begin{equation}
\label{cond1_GP}
	\alpha_{\infty}^h \geq \alpha_0 >0, \qqand \beta_{\infty}^h = 0,
\end{equation}
which are identical to the conditions given by Eq. \eqref{cond1}, \eqref{eq:diff_equation}.\\

Therefore Theorem~\ref{th:Linfty_G} follows as in Section~\ref{sec:L1G} because $ - \frac {1}{\widetilde \alpha_0} \frac{\bar \delta_\infty^h e^{\bar \delta_\infty^h t} }{e^{\bar \delta_\infty^h t}-1}$ is a sub-solution of the corresponding equation.\\

\bigskip

\section{$L^2$-type estimate}
\label{sec:L2}

We now investigate the $L^2$ space which has been used in situations where the $L^\infty$ estimate cannot be applied because the growth term depends on other quantities and cannot be differentiated with uniform control. As we shall see, the advantage of working in $L^2$ is to provide additional dissipation terms which do not appear in $L^1$ or $L^\infty$ while keeping an estimate compatible with the free boundary in opposition to $L
^1$. 

We proceed again by departing from Eq. \eqref{eq:Deltaw} which we write as
\begin{align}
\label{eq:estimateL2_vG}
\partial_t \frac{\abs{w}_{-}^{2}}{2} \leq -\frac{2}{N} \abs{w}_{-}^3 +   \nabla p\cdot \nabla  {\abs{w}_{-}^2} + \abs{w}_{-} \Delta \prt*{q \abs{w}_{-}} + \left[G'q- \frac 4 N G\right] \abs{w}_{-}^{2} .
\end{align}
And we distinguish the two cases $G\equiv 0$ or not.

\subsection{$L^2$-estimates when $G\equiv 0$. No weight.}
\label{subsec:L2_G0_noh}
$\mbox{}$\\

We are going to prove the following theorem
\begin{theorem} \label{th:L2G=0}
Assume $\inf_{0\leq p\leq p_M} \left( \frac{2}{N}-1 +\frac{q'}{2}\right) =: \alpha_0 \geq 0$ and  $q''\leq 0$, then
\[
\int_{\R^N} \abs{w(t)}_{-}^{2} \, \dx{x}  +2 \int_0^t \int_{\R^N} \prt*{\frac{2}{N}-1 +\frac{q'}{2}}\abs{w(s)}_{-}^{3}\, \dx{x} \dx{s} \leq \int_{\R^N} \abs{w(t=0)}_{-}^{2}  \, \dx{x} .
\]
When, for $t\in [0, T]$, we have $\alpha_0>0$, for solutions with compact support in $x$ , it holds
\[
\int_{\R^N} \abs{w(t)}_{-}^{2} \, \dx{x} \leq \frac {C(T)}{t^2} \qquad \forall t\in [0, T].
\]
\end{theorem}
For the Power-law, $q'= \gamma$, we recover the condition obtained in \cite{GPS19}, that is
\[
\gamma \geq 2-\frac 4 N. 
\]
In the case of the DHV pressure, \cf Eq. \eqref{eq:SingularPressure}, the  assumption $q'' \leq 0$ is not met.

Notice also that the argument can be localized, see \cite{david2020free}, and this allows to remove the compact support assumption.

\begin{proof}
Integrating Eq. \eqref{eq:estimateL2_vG}, and using two integration by parts, we obtain
\begin{align}
    \label{eq:TowardsL2Estimate_Neg_v}
	\ddt \int_{\R^N}  \frac{\abs{w}_{-}^{2}}{2} \, \dx{x} \leq -\frac{2}{N} \int_{\R^N} \abs{w}_{-}^3 \dx{x}  + \int_{\R^N} \abs{w}_{-}^3 \dx{x} 
	\underbrace{- \int_{\R^N} \nabla \abs{w}_{-}\nabla(q \abs{w}_{-}) \, \dx{x}}_{I} .
\end{align}
Integrating by parts, the first term, $I$, can be rewritten, \ie, 
\begin{align*}
	I & = - \int_{\R^N} \abs{w}_{-} \nabla q  \cdot \nabla \abs{w}_{-} \, \dx{x} - \int_{\R^N} q\abs{\nabla \abs{w}_{-}}^2\,\dx{x} \\
	&=  \frac{1}{2}\int_{\R^N}\Delta q \; \abs{w}_{-}^2 \, \dx{x} - \int_{\R^N} q\abs{\nabla \abs{w}_{-}}^2\,\dx{x}\\
	&= \frac{1}{2}\int_{\R^N} q'' \abs{\nabla p}^2 \abs{w}_{-}^2 \, \dx{x} -\frac{1}{2} \int_{\R^N} q' \abs{w}_{-}^3 \, \dx{x}  - \int_{\R^N} q\abs{\nabla \abs{w}_{-}}^2\,\dx{x},
\end{align*}
where we used the chain rule $\Delta q = q'' \abs{\nabla p}^{2}+q' w$.

Substituting the simplified expressions for $I$ into Eq. \eqref{eq:TowardsL2Estimate_Neg_v}, we obtain
\begin{align*}
    \ddt \int_{\R^N} \frac{\abs{w}_{-}^{2}}{2} \, \dx{x}  \leq - \int_{\R^N} \prt*{\frac{2}{N}-1 +\frac{q'}{2}}\abs{w}_{-}^{3}\, \dx{x} + \frac{1}{2} \int_{\R^N} q''\abs{\nabla p}^2 \abs{w}_{-}^{2}\, \dx{x}- \int_{\R^N} q\abs{\nabla \abs{w}_{-}}^2\,\dx{x}.
\end{align*}
Since $q\geq 0$, this ensures that the estimate can be closed if $q' > 2 - \frac 4 N $, for all $p \in (0,p_M)$ and  $q'' \leq 0$. Indeed, we deduce the inequality 
\begin{equation} \label{eq:L2NoNo}
\int_{\R^N} \abs{w(t)}_{-}^{2} \, \dx{x} \leq -  \int_{\R^N} \alpha_0  \abs{w(t)}_{-}^{3}.
\end{equation}
The conclusions of the theorem follow by time integration (first estimate) or using a sub-solution $\frac{C(T)}{t^2}$ (second estimate with regularizing effect).

\end{proof}

It is interesting to investigate if adding a weight can help us to include more general pressure laws.
\vspace{1em}

\subsection{$L^2$-estimates when $G\equiv 0$. With weights.}
\label{sec:L2estimatev}
$\mbox{}$\\
In order to be as general as possible, we add a weight $h = h(p)$.
\begin{theorem}
\label{th:G0_h}
Assume there exists a positive weight, $0<c \leq h(p)\leq c^{-1}$ such that the two differential inequalities
\begin{align*}
    \alpha_2^h := \frac{4}{N}+\frac{2 h'q}{h}-2 +q' \geq 0, \qqand \beta_2^h := \frac{h''q+q''h-h'}{h} \leq 0.
\end{align*}
are met. Then
\[
    \int_{\R^N} h|w|_-^2(t) \, \dx{x} +\int_0^t \int_{\R^N} \alpha_2^h h\wm^3\, \dx{x} \dx{t}\leq\int_{\R^N} h|w|_-^2(0) \, \dx{x}.
\]
\end{theorem}
Note that the same regularisation effect as in Theorem \ref{th:L2G=0} can be obtained if $\alpha_2^h> \alpha_0$, for some constant $\alpha_0>0$.
\begin{proof}
 Using the equation satisfied by the pressure, Eq. \eqref{eq:p_growth}, and Eq. \eqref{eq:estimateL2_vG} for $|w|_-^2$, we compute the evolution in time of $h|w|_-^2$ 
\begin{align}
\label{eq:evolution_hw_L2}
    \begin{split}
    \fpartial t \prt*{h|w|_-^2}
    &= h'|w|_-^2 \prt*{\abs{\nabla p}^2+ qw} + h  \partial_t |w|_-^2\\
    &\leq h' \abs{\nabla p}^2 |w|_-^2 - h' q |w|_-^3 \\
    &\qquad +2 h \prt*{-\frac{2}{N} |w|_-^3 + 2 |w|_- \nabla p \cdot \nabla |w|_- + |w|_- \Delta \prt*{q |w|_-}}.
    \end{split}
\end{align}
Integrating over $\R^{N}$ and with an integration by parts for the last two terms, we get
\begin{align}
\label{eq:evolutionL2_v}
\begin{split}
    \ddt \int_{\R^N} \prt*{h|w|_-^2} \, \dx{x} &
    \leq \int_{\R^N} h'\abs{\nabla p}^2 \wm^2 \, \dx{x} - \int_{\R^N} \prt*{\frac{4}{N}+\frac{h'q}{h}} h\wm^3\, \dx{x}\\
    &\quad +4 \int_{\R^N} h \wm \nabla p\cdot \nabla \wm \, \dx{x} - 2 \int_{\R^N} \nabla \prt*{h \wm} \nabla \prt*{q\wm}\, \dx{x} 
    \\
    &\leq -\int_{\R^N} h'\abs{\nabla p}^2 \wm^2 \dx{x} - \int_{\R^N}  \prt*{\frac{4}{N}- 2+\frac{h'q}{h}} h\wm^3  \dx{x}
    \\
    & \quad - 2\underbrace{\int_{\R^N} \nabla \prt*{h \wm} \nabla \prt*{q\wm}\, \dx{x}}_{I}.
\end{split}
\end{align}
Next, we need to address the term $I$. We compute
\begin{align} 
\nonumber
    I
    &= \int_{\R^N} \big[h' \nabla p \;\wm +h \nabla \wm \big]\cdot \big[q' \nabla p \;\wm +q \nabla \wm \big]  \dx{x}\\\label{eq:I_L2}
    &= \int_{\R^N} \left[ h'q' |\nabla p|^2 \;\wm^2 + \frac{h'q+hq'}{2}  \nabla p\cdot\nabla \wm^2 +hq |\nabla \wm|^2\right] \dx{x} \\\nonumber
    &= \int_{\R^N} \left[ h'q' |\nabla p|^2 \;\wm^2 + \frac{h'q+hq'}{2} \wm^3
    - \left[\frac{h'q+hq'}{2}\right]'  |\nabla p|^2 \wm^2
    +hq |\nabla \wm|^2\right] \dx{x} .
\end{align}
Reorganising the terms, we get 
\begin{align*}
I = \int_{\R^N} \left[\frac{h'q+hq'}{2} \wm^3 - \frac{h''q+hq''}{2}  |\nabla p|^2 \wm^2
 +hq |\nabla \wm|^2 \right] \dx{x}  
\end{align*}

Finally, substituting $I$ into Eq. \eqref{eq:evolutionL2_v}, we obtain
\begin{align*}
    \ddt\int_{\R^N}\prt*{h|w|_-^2} \, \dx{x} &\leq - \int_{\R^N} \alpha_2^h h\wm^3\, \dx{x} + \int_{\R^N} \beta_2^h \abs{\nabla p}^2 h\wm^2\, \dx{x}-2\int_{\R^N}  hq |\nabla \wm|^2 \dx{x}  ,
\end{align*}
where
\begin{align*}
    \alpha_2^h = \frac{4}{N}+\frac{2 h'q}{h}-2 +q', \qqand \beta_2^h = \frac{h''q+q''h-h'}{h}.
\end{align*}
By assumption 
\begin{equation}
\label{eq:condL2_v}
    \alpha_{2}^h \geq  0, \qqand \beta_{2}^h \leq 0.
\end{equation}
the statement holds true.
\end{proof}
\begin{remark}
\label{rem:L2_v_PL}
For the condition $\alpha_{2}^h\geq  0 $, the weight $h$ does not help when it is positive because the term generated by the weight, $\frac{2 h'q}{h^3}$ vanishes for $p=0$. However, loosing on the estimate near the free boundary, for instance  we may choose $h=p$, then we improve the range of possible parameters.  For  Power-law, we reach the conditions
\[
 \alpha_{2}^h = \frac{4}{p^2 N}+\frac{3\gamma}{p^2}-\frac{2}{p^2} \geq 0,  
\]
which is less restrictive than when $h=1$, while the condition $\beta_2^h \leq 0$ is fulfilled.

It is unclear to us how to choose the weight $h$ for DHV law in $L^2$.
\end{remark}

\vspace{1em}

Adding the reaction term $G$, we do not gain anything. Precisely, without the weight $h$, we obtain the same conditions on $q'$ and $q''$ as the ones in Theorem \ref{th:L2G=0}.  Since $G$ is bounded and decreasing, \cf Eq. \eqref{eq:AssumptionsOnG}, the term which involves $G$ can always be controlled. We decide not to report all the calculations because they can be derived easily from the ones in Section \ref{subsec:L2_G0_noh}.

\vspace{0.5cm}

To conclude the paper, we decide to add the reaction term $G>0$ and we can derive the $L^2$-estimate for $h\wm$.

\vspace{1em}

\subsection{$L^2$-estimates when $G\neq0$. With weights}
$\mbox{}$\\
This section is dedicated to proving the following theorem.
\begin{theorem}
\label{th:L2_Gno0}
Assume there exists a positive weight, $0<c \leq h(p)\leq c^{-1}$ such that the two differential inequalities  $\alpha_2^h  \geq 0, \beta_2^h  \leq 0$, are satisfied, with $\alpha_2^h$ and $\beta_2^h$ defined as in Theorem \ref{th:G0_h}.

Then there holds
\[
    \int_{\R^N} h|w|_-^2(t) \, \dx{x} +\int_0^t \int_{\R^N} \alpha_2^h h\wm^3\, \dx{x} \dx{t}\leq\int_{\R^N} h|w|_-^2(0) \, \dx{x}  +\int_0^t\int_{\R^N} \bar \delta_2^h h\wm^2 \, \dx{x}\dx{t},
\]
where 
\begin{align*}
    &\bar \delta_2^h = \sup_{0\leq p \leq  p_M}\left\{2G' q+G\prt*{2\left(1-\frac{4}{N}\right) - \frac{h'q+q'h}{h}}\right\}.
\end{align*}
\end{theorem}
As before, the same regularisation effect as in Theorem \ref{th:L2G=0} is obtained if we can guarantee that $\alpha_2^h> \alpha_0$, for a positive constant $\alpha_0>0$. Theorem \ref{th:L2_Gno0} proves that adding the reaction, the estimate does not gain anything.

\begin{proof}
Starting from Eq. \eqref{eq:evolution_hw_L2} and including the additional growth terms $G$, the evolution of $h\wm^2$ becomes
\begin{align}
\label{eq:evolution_hw_L2_G}
    \begin{split}
    \fpartial t \prt*{h|w|_-^2}
    &\leq h' \abs{\nabla p}^2 |w|_-^2 - h' q |w|_-^3 \\
    &\quad +2 h \prt*{-\frac{2}{N} |w|_-^3 + \nabla p \cdot \nabla \wm^2 + |w|_- \Delta \prt*{q |w|_-}}\\
    &\quad+ 2 h\left(G' q -\frac{4}{N} G\right)\wm^2.
    \end{split}
\end{align}
Integrating in space and by an integration by parts, we get
\begin{align}
\label{eq:L2_of_hw_neg_sqd}
\begin{split}
	\ddt \int_{\R^N} h\wm^2 \, \dx{x} 
	&\leq -\int_{\R^N} \left(\frac{4}{ N} +\frac{ qh'}{h} -2\right)h\wm^3\, \dx{x} - \int_{\R^N} h' \wm^2 \abs{\nabla p}^2 \, \dx{x}\\
	&\quad  - 2\underbrace{\int_{\R^N} \nabla\left(h\wm\right)\cdot\nabla\left(q\wm\right)\, \dx{x}}_{{I}} + 2\int_{\R^N} h\wm^2 \left(G' q + G\left(1-\frac{4}{N}\right)\right) \dx{x}.
\end{split}
\end{align}

Using Eq. \eqref{eq:I_L2} and the definition of $w$, \cf Eq. \eqref{eq:DefOf_w}, the term ${I}$ simplifies to
\begin{align}
\label{eq:termC_hw_L2}
\begin{split}
	{I} 	&= \frac{1}{2}\int_{\R^N}  \left(h' q + h q'\right)\wm^3 \, \dx{x} - \frac{1}{2}\int_{\R^N}\left(h'' q + h q''\right) \abs{\nabla p}^2 \wm^2\, \dx{x}\\
	&\quad+ \frac{1}{2}\int_{\R^N} 
	\left(h' q + h q'\right)G \wm^2\, \dx{x} +\int_{\R^N} hq \abs{\nabla \wm}^2\, \dx{x}
\end{split}	
\end{align}

Substituting everything into Eq. \eqref{eq:L2_of_hw_neg_sqd},  we obtain
\begin{align*}
		\ddt \int_{\R^N} \frac{h\wm^2}{2} \, \dx{x} &\leq -\int_{\R^N} \alpha_2^hh\wm^3\, \dx{x} + \int_{\R^N} \abs{\nabla p}^2 \beta_2^hh\wm^2 \, \dx{x} + \int_{\R^N} \delta_2^h h\wm^2 \, \dx{x},
\end{align*}
where
\begin{align*}
    \alpha_2^h = \frac{4}{N}+\frac{2 qh'}{h}
    &-2+q', \qqand \beta_2^h = \frac{h''q+q''h-h'}{h},
\end{align*}
as well as
\begin{align*}
    &\delta_2^h = 2G' q+G\prt*{2\left(1-\frac{4}{N}\right) - \frac{h'q+q'h}{h}} \leq \bar{\delta}_2<\infty,
\end{align*}
where $\bar{\delta}_2$ is as in the  statement.
By assumption
\begin{equation}
    \label{eq:systemhw_L2}
    \alpha_2^h \geq 0 \qqand \beta_2^h \leq 0,
\end{equation}
and the statement holds true.
\end{proof}
\begin{remark}
 Different techniques have been used to perform an $L^2$-bound. For instance, in \cite{david2020free}, the authors compute exactly the square of the expression $(w-G)$ and then they apply the Young's inequality to the term
    $$
    2\left(1-\frac{2}{N}\right)G\wm^2.
    $$
    This choice allows them to get an additional small constant $\mu>0$ in the expression of $\alpha_2^h$ which can help in getting a weaken condition to close the estimate. However, there is an additional term which involves the growth term, \ie,
    $$
    C \int_{\R^N} G^2 \wm\, \dx{x},
    $$
    with $C>0$ a constant. Now the estimate can only be closed by assuming the standard conditions on $G$, \cf Eq. \eqref{eq:AssumptionsOnG}, plus an additional condition on the domain. Finally, the expression of $\delta_2^h$ change a bit, \ie,
    $$
    \left(\delta_2^h\right)_Y = \frac{2G'q}{h} - G \left(\frac{h'q+q'h}{h^2}\right) \leq \Tilde{\delta}_2 < \infty,
    $$
    which can be bounded by controlling $h'$.
\end{remark}

\begin{remark}
     In the Power-law case, regarding the conditions on $\alpha_2^h$ and $\beta_2^h$, we can refer to Remark \ref{rem:L2_v_PL}. The last one, substituting $ q = p \gamma$ and $h(p) = c_1 +p$, $\delta_2^h$ becomes
     \begin{align}
         \delta_2^h = \gamma p \prt*{2G'-c_1 G} - \prt{p+c_1} \prt*{G \gamma + 2 G} < \infty.
     \end{align}
     Note that this expression is always non-positive for non-negative and decreasing growth terms, \cf Eq. \eqref{eq:AssumptionsOnG}.
\end{remark}

\bigskip

\section{Conclusions}
\label{sec:conclusions}

The Aronson-B\'enilan estimate has proven to be a fundamental tool in order to study regularity and asymptotic in several problems related to the porous media equations. Even if it has been used mainly to control the Laplacian of the pressure from below by a term as $\Delta p(t) \geq -\frac C t $, one may use it in other Lebesgue spaces. We have systematically studied the restrictions on the parameters and the conclusions that one can draw in $L
^1$, $L^2$ and $L^\infty$ (original work of Aronson and B\'enilan). In particular we considered two specific forms of the pressure law, the Power-law and DHV law.

Our conclusions are that the $L^\infty$ setting provides the widest range of parameters, generating the strongest estimate. For instance, it can be applied to both pressure laws. 
On the other hand, the $L^2$ estimate requires restrictions on the parameters (which exclude DHV law) but is enough to estimate the Laplacian of the pressure for the Power-law. Because of integration by parts, and because a dissipation term occurs explicitly, it is however useful for some strongly coupled problem where $L^\infty$ bounds are not possible. The $L^1$ estimate turns out to be the simplest but is only useful in space dimension $N=1$.

When weights are included in order to treat more general equations of state, we improve the results in \cite{CP82} and we obtain estimates correctly scaled with respect to the Hele-Shaw limit, which, with our notations, is expressed as $|w|_- \approx 0$ for $\gamma \gg 1$ or $\epsilon \ll 1$.

If one wishes to estimate the quantity $p \Delta p(t)$, loosing regularity near the free boundary, then one can  drastically extend the range of possible pressure laws.

\section*{Acknowledgements}
Part of this work was done while G.B. was visiting the Laboratoire Jacques-Louis Lions at the Sorbonne Universit\'e whose hospitality is gratefully acknowledged. The research experience was supported by LIA-LYSM: AMU-CNRS.ECM-INdAM funding. The work of G.B. was also partially supported by GNFM-INdAM.

B.P. has received funding from the European Research Council (ERC) under the European Union's Horizon 2020 research and innovation programme (grant agreement No 740623).

M.S. fondly acknowledges the support of the Fondation Sciences Math\'ematiques de Paris (FSMP) for the postdoctoral fellowship.

\bibliographystyle{plain}

\begin{thebibliography}{10}
\bibitem{A99_brinkman}
G.\ Allaire.
\newblock Homogenization of the navier-stokes equations and derivation of
  brinkman's law.
\newblock {\em Math{\'e}matiques appliqu{\'e}es aux sciences de
  l'ing{\'e}nieur (Santiago, 1989)}, pages 7--20, 1991.



\bibitem{AA95}
S.B. Angenent and D.G. Aronson.
\newblock {The focusing problem for the radially symmetric porous medium
  equation}.
\newblock {\em Communications in Partial Differential Equations},
  20(7-8):1217--1240, 1995.
  
\bibitem{Aro69}
D.~G. Aronson.
\newblock Regularity propeties of flows through porous media.
\newblock {\em SIAM Journal on Applied Mathematics}, 17(2):461--467, 1969.

\bibitem{Aro70}
D.~G. Aronson.
\newblock Regularity properties of flows through porous media: A counterexample.
\newblock 19(2):299-307, 1970.

\bibitem{Aro70a}
D.~G. Aronson.
\newblock Regularity properties of flows through porous media: The interface.
\newblock 37(1):1-10, 1970.


\bibitem{AB79}
D.~G. Aronson and P.~B{\'e}nilan.
\newblock R{\'e}gularit{\'e} des solutions de l'{\'e}quation des milieux poreux
  dans $\mathbb{R}^n$.
\newblock {\em CR Acad. Sci. Paris S{\'e}r. AB}, 288(2):A103--A105, 1979.


\bibitem{AG93}
D.~G. Aronson and J.~Graveleau.
\newblock {A selfsimilar solution to the focusing problem for the porous medium
  equation}.
\newblock {\em European Journal of Applied Mathematics}, 4(1):65--81, 1993.



\bibitem{Bar52}
G.~I. Barenblatt.
\newblock On some unsteady motions of a liquid and gas in a porous medium.
\newblock {\em Akad. Nauk SSSR. Prikl. Mat. Meh.}, 16:67--78, 1952.



\bibitem{bellomo2009complexity}
N.~Bellomo, C.~Bianca, and M.~Delitala.
\newblock Complexity analysis and mathematical tools towards the modelling of
  living systems.
\newblock {\em Physics of Life Reviews}, 6(3):144--175, 2009.

\bibitem{bellomo2003modeling}
N.~Bellomo and E.~De~Agelis.
\newblock Modeling and simulation of tumor development, treatment, and control.
\newblock {\em Mathematical and computer modelling}, 37(11), 2003.

\bibitem{bellomo2008foundations}
N.~Bellomo, N.K. Li, and P.~K. Maini.
\newblock On the foundations of cancer modelling: selected topics,
  speculations, and perspectives.
\newblock {\em Mathematical Models and Methods in Applied Sciences},
  18(04):593--646, 2008.

\bibitem{benilan1989limit}
P.~B{\'e}nilan, L.~Boccardo, and M.~A. Herrero.
\newblock On the limit of solutions of $u_t= \delta u^{m}$ as $m \to \infty$.
\newblock 1989.

\bibitem{BCGRS10}
D.~Bresch, T.~Colin, E.~Grenier, B.~Ribba, and O.~Saut.
\newblock Computational modeling of solid tumor growth: the avascular stage.
\newblock {\em SIAM Journal on Scientific Computing}, 32(4):2321--2344, 2010.




\bibitem{BPPS19}
F.~Bubba, B.~Perthame, C.~Pouchol, and M.~Schmidtchen.
\newblock Hele--shaw limit for a system of two reaction-(cross-) diffusion
  equations for living tissues.
\newblock {\em Archive for Rational Mechanics and Analysis}, 236(2):735--766,
  2020.

\bibitem{BD09}
H.~Byrne and D.~Drasdo.
\newblock Individual-based and continuum models of growing cell populations: a
  comparison.
\newblock {\em Journal of mathematical biology}, 58(4-5):657, 2009.

\bibitem{byrne1996modelling}
H.~M. Byrne and M.~A.~J. Chaplain.
\newblock Modelling the role of cell-cell adhesion in the growth and
  development of carcinomas.
\newblock {\em Mathematical and Computer Modelling}, 24(12):1--17, 1996.

\bibitem{CE83}
L.A. Caffarelli and L.C. Evans.
\newblock Continuity of the temperature in the two-phase stefan problem.
\newblock {\em Archive for Rational Mechanics and Analysis}, 81(3):199--220,
  1983.

\bibitem{CF79}
L.A. Caffarelli and A.~Friedman.
\newblock {Continuity of the Density of a Gas Flow in a Porous Medium}.
\newblock {\em Transactions of the American Mathematical Society}, 252:99,
  1979.

\bibitem{CF79a}
L.A. Caffarelli and A.~Friedman.
\newblock Regularity of the free boundary for the one-dimensional flow of gas
  in a porous medium.
\newblock {\em American Journal of Mathematics}, 101(6):1193--1218, 1979.

\bibitem{CF80}
L.A. Caffarelli and A.~Friedman.
\newblock {Regularity of the Free Boundary of a Gas Flow in an n-dimensional
  Porous Medium}.
\newblock {\em Indiana University Mathematics Journal}, 29:361--391, 1980.




\bibitem{CVW87}
L.A. Caffarelli, J.-L. V{\'a}zquez, and N.I. Wolanski.
\newblock {Lipschitz continuity of solutions and interfaces of the
  n--dimensional porous medium equation}.
\newblock {\em Indiana University mathematics journal}, 36(2):373--401, 1987.

\bibitem{CFSS17}
J.~A. Carrillo, S.~Fagioli, F.~Santambrogio, and M.~Schmidtchen.
\newblock Splitting schemes and segregation in reaction cross-diffusion
  systems.
\newblock {\em SIAM Journal on Mathematical Analysis}, 50(5):5695--5718, 2018.

\bibitem{CDHV19}
A.~Chertock, P.~Degond, S.~Hecht, and J.P. Vincent.
\newblock Incompressible limit of a continuum model of tissue growth with
  segregation for two cell populations.
\newblock {\em Mathematical biosciences and engineering: MBE}, 16(5):5804,
  2019.

\bibitem{ciarletta2011radial}
P.~Ciarletta, L.~Foret, and M.~Ben~Amar.
\newblock The radial growth phase of malignant melanoma: multi-phase modelling,
  numerical simulations and linear stability analysis.
\newblock {\em Journal of the Royal Society Interface}, 8(56):345--368, 2011.

\bibitem{CP82}
M.~G. Crandall and M.~Pierre.
\newblock Regularizing effects for {$u_{t}=\Delta \varphi (u)$}.
\newblock {\em Trans. Amer. Math. Soc.}, 274(1):159--168, 1982.

\bibitem{cui2008asymptotic}
S.~Cui and J.~Escher.
\newblock Asymptotic behaviour of solutions of a multidimensional moving
  boundary problem modeling tumor growth.
\newblock {\em Communications in Partial Differential Equations},
  33(4):636--655, 2008.
  
  
  
\bibitem{Dar1856}
H.~P.~G. Darcy.
\newblock {\em Les Fontaines publiques de la ville de Dijon. Exposition et
  application des principes {\`a} suivre et des formules {\`a} employer dans
  les questions de distribution d'eau, etc}.
\newblock V. Dalamont, 1856.

\bibitem{david2020free}
N.~David and B.~Perthame.
\newblock Free boundary limit of tumor growth model with nutrient.
\newblock {\em arXiv preprint arXiv:2003.10731}, 2020.

\bibitem{DPSV20}
T.~D{\k{e}}biec, B.~Perthame, M.~Schmidtchen, and N.~Vauchelet.
\newblock Incompressible limit for a two-species model with coupling through
  brinkman's law in any dimension.
\newblock {\em ArXiv}, 2020.



\bibitem{DS20}
T.~D{\k{e}}biec and M.~Schmidtchen.
\newblock Incompressible limit for a two-species tumour model with coupling
  through brinkman's law in one dimension.
\newblock {\em Acta Applicandae Mathematicae}, pages 1--19, 2020.


\bibitem{DHV20}
P.~Degond, S.~Hecht, N.~Vauchelet, and and.
\newblock Incompressible limit of a continuum model of tissue growth for two
  cell populations.
\newblock {\em Networks {\&} Heterogeneous Media}, 15(1):57--85, 2020.


\bibitem{deroulers2009modeling}
C. Deroulers, M. Aubert, M. Badoual, and B. Grammaticos.
\newblock Modeling tumor cell migration: from microscopic to macroscopic
  models.
\newblock {\em Physical Review E}, 79(3):031917, 2009.


\bibitem{DJ20}
P.-E. Druet and A.~J\"ungel.
\newblock Analysis of cross-diffusion systems for fluid mixtures driven by a
  pressure gradient.
\newblock {\em SIAM Journal on Mathematical Analysis}, 52(2):2179--2197, 2020.


\bibitem{friedman2004hierarchy}
A.~Friedman.
\newblock A hierarchy of cancer models and their mathematical challenges.
\newblock {\em Discrete \& Continuous Dynamical Systems-B}, 4(1):147, 2004.

\bibitem{friedman2008stability}
A.~Friedman and B.~Hu.
\newblock Stability and instability of {Liapunov-Schmidt} and {Hopf}
  bifurcation for a free boundary problem arising in a tumor model.
\newblock {\em Transactions of the American Mathematical Society},
  360(10):5291--5342, 2008.

\bibitem{GST2019}
B.~Gess, J.~Sauer, and E.~Tadmor.
\newblock Optimal regularity in time and space for the porous medium equation.
\newblock {\em arXiv preprint arXiv:1902.08632}, 2019.

\bibitem{G2017}
B.~Gess.
\newblock Optimal regularity for the porous medium equation.
\newblock {\em arXiv preprint arXiv:1708.04408}, 2017.



\bibitem{gil2001convergence}
O.~Gil and F.~Quir{\'o}s.
\newblock Convergence of the porous media equation to {Hele-Shaw}.
\newblock {\em Nonlinear Analysis: Theory, Methods \& Applications},
  44(8):1111--1131, 2001.

\bibitem{gil2003boundary}
O.~Gil and F.~Quir{\'o}s.
\newblock Boundary layer formation in the transition from the porous media
  equation to a {Hele-Shaw} flow.
\newblock In {\em Annales de l'IHP Analyse non lin{\'e}aire}, volume~20, pages
  13--36, 2003.

\bibitem{greenspan1972models}
H.~P. Greenspan.
\newblock Models for the growth of a solid tumor by diffusion.
\newblock {\em Studies in Applied Mathematics}, 51(4):317--340, 1972.

\bibitem{greenspan1976growth}
H.~P. Greenspan.
\newblock On the growth and stability of cell cultures and solid tumors.
\newblock {\em Journal of theoretical biology}, 56(1):229--242, 1976.

\bibitem{GN75}
W.~S.~C. Gurney and R.~M. Nisbet.
\newblock The regulation of inhomogeneous populations.
\newblock {\em Journal of Theoretical Biology}, 52(2):441--457, 1975.



\bibitem{GPS19}
P.~Gwiazda, B.~Perthame, and A.~{\'S}wierczewska-Gwiazda.
\newblock A two-species hyperbolic--parabolic model of tissue growth.
\newblock {\em Communications in Partial Differential Equations},
  44(12):1605--1618, 2019.

\bibitem{HV17}
S.~Hecht and N.~Vauchelet.
\newblock {Incompressible limit of a mechanical model for tissue growth with
  non-overlapping constraint}.
\newblock {\em Communications in Mathematical Sciences}, 15(7):1913--1932,
  2017.



\bibitem{IKO62}
A.~M. Il'in, A.~S. Kalashnikov, and O.~A. Oleinik.
\newblock Linear equations of the second order of parabolic type.
\newblock 17(3):1-143, 1962.


\bibitem{Ito52}
Y.~It{\^o}.
\newblock The growth form of populations in some aphids, with special reference
  to the relation between population density and the movements.
\newblock {\em Researches on Population Ecology}, 1(1):36--48, 1952.




\bibitem{J10}
P.~E. Jabin.
\newblock Differential equations with singular fields.
\newblock {\em Journal de math{\'e}matiques pures et appliqu{\'e}es},
  94(6):597--621, 2010.

\bibitem{Kal66}
A.~S. Kalashnikov.
\newblock The occurrence of singularities in solutions of the non-steady
  seepage equation.
\newblock {\em {USSR} Computational Mathematics and Mathematical Physics},
  7(2):269--275, jan 1967.

\bibitem{KP15}
Inwon Kim and Norbert Pozar.
\newblock {Porous medium equation to Hele-Shaw flow with general initial
  density}.
\newblock 2015.

\bibitem{LTWZ19}
J.-G. Liu, M.~Tang, L.~Wang, and Z.~Zhou.
\newblock Towards understanding the boundary propagation speeds in tumor growth
  models.
\newblock {\em arXiv preprint arXiv:1910.11502}, 2019.

\bibitem{lowengrub2009nonlinear}
J.~S. Lowengrub, H.~B. Frieboes, F.~Jin, Y.-L. Chuang, X.~Li, P.~Macklin, S.~M.
  Wise, and V.~Cristini.
\newblock Nonlinear modelling of cancer: bridging the gap between cells and
  tumours.
\newblock {\em Nonlinearity}, 23(1):R1, 2009.

\bibitem{MPQ17}
A.~Mellet, B.~Perthame, and F.~Quir\'{o}s.
\newblock A {H}ele-{S}haw problem for tumor growth.
\newblock {\em J. Funct. Anal.}, 273(10):3061--3093, 2017.

\bibitem{Mor50}
M.~Morisita.
\newblock Population density and dispersal of a water strider. gerris
  lacustris: Observations and considerations on animal aggregations.
\newblock {\em Contributions on Physiology and Ecology, Kyoto University},
  65:1--149, 1950.




\bibitem{OKYL58}
O.~A. Oleinik, A.~S. Kalasinkov, and Y.L. Czou.
\newblock The {C}auchy problem and boundary problems for equations of the type
  of non-stationary filtration.
\newblock {\em Izv. Akad. Nauk SSSR. Ser. Mat.}, 22:667--704, 1958.

\bibitem{Pat59}
R.~E. Pattle.
\newblock Diffusion from an instantaneous point source with a
  concentration-dependent coefficient.
\newblock {\em The Quarterly Journal of Mechanics and Applied Mathematics},
  12(4):407--409, 1959.

\bibitem{PQTV14}
B.~Perthame, F.~Quir{\'o}s, M.~Tang, and N.~Vauchelet.
\newblock Derivation of a hele-shaw type system from a cell model with active
  motion.
\newblock {\em Interfaces and Free Boundaries}, 16:489--508, 2014.

\bibitem{PQV14}
B.~Perthame, F.~Quir\'{o}s, and J.-L. V\'{a}zquez.
\newblock The {H}ele-{S}haw asymptotics for mechanical models of tumor growth.
\newblock {\em Archive for Rational Mechanics and Analysis}, 212(1):93--127,
  2014.

\bibitem{PV15}
B.~Perthame and N.~Vauchelet.
\newblock Incompressible limit of a mechanical model of tumour growth with
  viscosity.
\newblock {\em Phil. Trans. R. Soc. A}, 373(2050):20140283, 2015.

\bibitem{preziosi2009multiphase}
L.~Preziosi and A.~Tosin.
\newblock Multiphase modelling of tumour growth and extracellular matrix
  interaction: mathematical tools and applications.
\newblock {\em Journal of mathematical biology}, 58(4-5):625, 2009.

\bibitem{PX20}
B.C. Price and X.~Xu.
\newblock Global existence theorem for a model governing the motion of two cell
  populations.
\newblock {\em arXiv preprint arXiv:2004.05939}, 2020.

\bibitem{RBEJPJ10}
J.~Ranft, M.~Basana, J.~Elgeti, J.F. Joanny, J.~Prost, and F.~J\"ulicher.
\newblock Fluidization of tissues by cell division and apoptosis.
\newblock {\em Natl. Acad. Sci. USA}, 49:657--687, 2010.




\bibitem{ribba2006multiscale}
B.~Ribba, O.~Saut, T.~Colin, D.~Bresch, E.~Grenier, and J.-P. Boissel.
\newblock A multiscale mathematical model of avascular tumor growth to
  investigate the therapeutic benefit of anti-invasive agents.
\newblock {\em Journal of theoretical biology}, 243(4):532--541, 2006.

\bibitem{Sab61}
E.~S. Sabinina.
\newblock On the cauchy problem for the equation of nonstationary gas
  filtration in several space variables.
\newblock In {\em Doklady Akademii Nauk}, volume 136, pages 1034--1037. Russian
  Academy of Sciences, 1961.

\bibitem{Vaz07}
J.-L. V{\'a}zquez.
\newblock {\em The porous medium equation: mathematical theory}.
\newblock Oxford University Press, 2007.

\bibitem{ZR66}
I.~B. Zel'dovich and I.~P. Raizer.
\newblock Physics of shock waves and high-temperature phenomena.
\newblock 1966.





\end{thebibliography}

\end{document}